\documentclass{amsart}

\usepackage{amsmath, amsfonts, amssymb, amscd}

\usepackage{amsthm}
\usepackage{amscd}
\newtheorem{thm}{Theorem}

\newtheorem{cor}{Corollary}
\newtheorem{lem}{Lemma}

\begin{document}

\title{Hilbert series of  algebras  associated to directed graphs} 

\author{Vladimir Retakh}
\address{Department of Mathematics, Rutgers University, Piscataway, NJ 08854-8019, USA}
\email{vretakh@math.rutgers.edu}
\author{Shirlei Serconek}
\address{IME-UFG, CX Postal 131, Goiania - GO, CEP 74001-970,  Brazil}
\email{serconek@math.rutgers.edu}
\author{Robert Lee Wilson}
\address{Department of Mathematics, Rutgers University, Piscataway, NJ 08854-8019, USA}
\email{rwilson@math.rutgers.edu}

\keywords{Hilbert series, directed graphs, quadratic algebras}
\subjclass{05E05; 15A15; 16W30}

\begin{abstract}
We compute the Hilbert series of some algebras associated to directed graphs
and related to factorizations of noncommutative polynomials.
\end{abstract}

\maketitle

\section{Introduction}

In \cite{GRSW} we introduced a new class of algebras $A(\Gamma )$
associated to layered directed graphs $\Gamma $. These algebras arose 
as generalizations of the algebras $Q_n$ (which are related to 
factorizations of noncommutative polynomials, see  \cite{GGRSW,GRW,SW}),
but the new class of algebras seems to be interesting by itself.

Various results have been proven for algebras $A(\Gamma )$.
In \cite{GRSW} we constructed a linear basis in $A(\Gamma )$.
In \cite{RSW} we showed that  algebras $A(\Gamma )$ are defined by
quadratic relations for a large class of directed graphs and proved
that in this case they are Koszul algebras. It follows immediately that the dual
algebras to $A(\Gamma )$ are also Koszul and that their Hilbert series
are related.

In this paper we continue to study algebras $A(\Gamma )$. 
In Section 2 we recall the definition of the algebra $A(\Gamma)$ and the 
construction of a basis for $A(\Gamma)$ given in \cite{GRSW}.  
In Section 3 we prove the main result of the paper, an expression for the 
Hilbert series, $H(A(\Gamma),t)$ of the algebra $A(\Gamma)$ corresponding 
to a layered graph $\Gamma$ with a unique element $*$ of level $0$.  
In stating this we denote the level of $v$ by $|v|$ and write 
$v > w$ to indicate that $v$ and $w$ are 
vertices of the directed graph $\Gamma$ and that there is a 
directed path from $v$ to $w$. Then we have:

$$H(A(\Gamma),t)  = \frac{1-t}{1 + \sum_{v_1 > v_2 \dots > v_{\ell}\geq *} 
(-1)^{\ell} t^{|v_1| - |v_{\ell}|+1}}.$$

The proof uses matrices $\zeta(t)$ and $\zeta(t)^{-1}$ which generalize 
the zeta function and the M\"obius function for partially ordered sets.

In Section 4 we specialize our results to the case of the 
Hasse graph of the lattice of subsets of  a finite set, giving a derivation  
of the Hilbert series for the algebras $Q_n$ that is shorter and more conceptual 
than that in \cite{GGRSW}. 
In Section 5 we treat the case of the Hasse graph of the lattice of subspaces 
of a finite-dimensional vector space over a
finite field. Finally, in Section 6, we define the complete layered graph 
$\mathbf{C}[m_n,m_{n-1},\dots ,m_1,m_0]$ and compute the Hilbert series of 
$A(\mathbf{C}[m_n,m_{n-1},\dots ,m_1,1])$.

During preparation of this paper Vladimir Retakh was partially supported 
by NSA.

\section{The algebra $A(\Gamma)$}

We begin by recalling  the definition of the algebra $A(\Gamma).$
Let  $\Gamma = (V, E)$ be a \textbf{directed graph}. That is, $V$ is a 
set (of vertices), 
$E$ is a set (of edges), and $\mathbf{ t}: E \rightarrow V$ and $\mathbf{ h}: E \rightarrow V$ are functions. ($\mathbf{ t}(e)$ is the {\it tail} 
of $e$ and $\mathbf{ h}(e)$ is the {\it head} of $e$.)       

We say that $\Gamma$ is \textbf{layered} if $V = \cup _{i=0}^n V_i$, $E = \cup_{i=1}^n E_i$, 
 $\mathbf{ t}: E_i \rightarrow V_i$, \ $\mathbf{ h}: E_i \rightarrow V_{i-1}$. 
If $v \in V_i$ we will write $|v| = i.$

We will assume throughout the remainder of the paper that $\Gamma = (V, E)$ is a layered graph with 
$V = \cup_{i=0}^n V_i$, that $V_0 = \{*\}$, and that, for every $v \in  V_+ = \cup_{i=1}^n V_i$, 
$\{ e \in E \ | \ \mathbf{ t}(e) = v \} \neq \emptyset$. For each $v \in V_+$ fix, arbitrarily, some $e_v \in E$ with ${\mathbf t}(e_v) = v$.  

If $v, \ w \in V$, a {\bf path} from $v$ to $w$ is a sequence of edges 
$\pi = \{ e_1, e_2, \dots,e_m \}$ with $\mathbf{ t}(e_{1}) = v$, 
$\mathbf{ 
h}(e_m) = w$ and $\mathbf{ t}(e_{i+1}) = \mathbf{ h}(e_i)$ for $1 \leq i < m$.  We write $v = \mathbf{ t}(\pi)$, $w = \mathbf{ h}(\pi)$.  
We also write $v > w$ if there  is a path from 
$v$ to $w$.  Define 
$P_{\pi}(\tau ) = (\tau -e_1)(\tau -e_2)\dots (\tau -e_m) \in 
T(E)[\tau ]$ and write
$$P_{\pi}(\tau ) = \sum_{j=0}^n e(\pi,j)\tau ^{m-j}.$$
Let $\pi_v$ denote the path $\{e_1,\dots ,e_{|v|}\}$ 
from $v$ to $*$ 
with $e_1 = e_{v}, e_{i+1} = e_{{\mathbf h}(e_i)}$ for $1 \le i < |v|$, and ${\mathbf h}(e_{|v|}) = *$.

Recall  that $R$ is the ideal of $T(E)$ generated by $$\{ e(\pi_1,k) - e(\pi_2,k) \ | \ \mathbf{ t}(\pi_1)=\mathbf{ t}(\pi_2), \mathbf{ h}(\pi_1)=\mathbf{ h}(\pi_2) , \ 1 \leq k \leq l(\pi_1) \}.$$
The algebra $A(\Gamma)$ is the quotient $T(E)/R$.  

For $v \in V_+$ and $1 \le k \le |v|$ we define $\hat{e}(v,k)$ to be the image in $A(\Gamma)$ of the product
$e_1\dots e_k$ in $T(E)$ where $\pi_v = \{e_1,\dots ,e_{|v|}\}$.  

If $(v,k), (u,l) \in V \times {\mathbf N}$ we say $(v,k)$  {\bf covers}  $(u,l)$ if $v > u$ and $k = |v| - 
|u|$.  In this case we write $(v,k) \gtrdot (u,l)$.  (In \cite{GRSW} we used different terminology and notation:  if $(v,l) \gtrdot (u,l)$ we said $(v,l) $ can be  composed with $(u,l)$ and wrote $(v,l) \models (u,l)$.)

The following theorem is proved in \cite [Corollary 4.5]{GRSW}.

\begin{thm}  Let $\Gamma = (V,E)$  be a layered graph, $V = \cup_{i=0}^n V_i, $ and $V_0 = \{*\}$ 
where $*$ is the unique minimal vertex of $\Gamma$. 
Then $$\{\hat{e}(v_1,k_1)\dots \hat{e}(v_{\ell},k_{\ell})|l \ge 0, 
v_1,\dots ,v_{\ell} \in 
V_+, 1 \le k_i \le |v_i|, (v_i,k_i) \not\gtrdot (v_{i+1},k_{i+1})\}$$
is a basis for $A(\Gamma).$
\end{thm}
\section{The Hilbert series of $A(\Gamma)$}

Let $h(t)$  denote the Hilbert series $H(A(\Gamma),t)$,  where $\Gamma$ is a layered  graph with unique minimal element $*$ of level $0$.
If $X \subseteq A(\Gamma)$ is a set of homogeneous elements 
(so  $X = \cup_{i=0}^{\infty} X_i$ where $X_i = X \cap A(\Gamma)_i$), 
denote the "graded cardinality" $\sum_{i=0}^{\infty} |X_i|t^i$ of $X$ by $||X||$.
Let $B$ denote the basis for $A(\Gamma)$ described in Theorem 1 
and, for $v \in V_+$, let $B_v = \{\hat{e}(v_1,k_1)\dots 
\hat{e}(v_{\ell},k_{\ell}) 
\in B|v_1 = v\}$.   Then $B = \{1\} \cup \bigcup_{v \in V_+} B_v.$   
Let $h_v(t)$ denote 
the graded dimension of the subspace of $A(\Gamma)$ spanned by $B_v$. 
Since $B$ is linearly independent, we have $||B|| = h(t)$ and $||B_v|| = h_v(t).$ 
Then
$$||B|| = h(t) = 1 + \sum_{v  \in V_+}h_v(t)$$  
Let $C_v = \bigcup_{k=1}^{|v|} \hat{e}(v,k)B$. Then 
$$||C_v|| = (t + \dots + t^{|v|})h(t) = t\left (\frac{t^{|v|} - 1}{t-1}\right 
)h(t).$$  
Now $C_v \supseteq B_v$.  Let $D_v$ denote the compliment of $B_v$ in $C_v$.  Then
$$D_v = \{\hat{e}(v,k)\hat{e}(v_1,k_1)\dots \hat{e}(v_{\ell},k_{ell})|1 
\le k \le 
|v|, $$ $$(v,k) \gtrdot (v_1,k_1), \hat{e}(v_1,k_1)\dots 
\hat{e}(v_{\ell},k_{\ell}) \in B\}  $$ 
and so 
$$D_v = \bigcup_{v > v_1 > *} \hat{e}(v,|v|-|v_1|)B_{v_1}.$$  Then $||D_v|| 
= \sum_{v > v_1 > *} t^{|v|-|v_1|}h_{v_1}(t)$
and so
$$h_v(t) = ||B_v|| = ||C_v|| - ||D_v|| = 
t\left (\frac{t^{|v|} - 1}{t-1}\right )h(t) - \sum_{v > w>*} t^{|v| - 
|w|}h_w(t).$$

This equation may be written in matrix form.  Arrange the elements of $V$ in nonincreasing
order and index the elements of vectors and matrices by this ordered set.  Let ${\mathbf h}(t)$
denote the column vector with entry $h_v(t)$ in the $v$-position (where we set $h_*(t)=1)$, let
$\mathbf u$ denote the vector with $t^{|v|}$ in the $v$-position,
$\mathbf e_*$ denote the vector with $\delta _{*v}$ in the $v$-position,
let ${\mathbf 1}$ denote the column vector all of whose entries are $1$,
and let $\zeta(t)$ denote the matrix with entries $\zeta_{v,w}(t)$ for $ v,w \in V$ where
$\zeta_{v,w}(t) = t^{|v|-|w|}$ if $v \ge w$ and $0$ otherwise.
Note that $$\zeta (t)\mathbf e_*=\mathbf u.$$

Then we have
$$\zeta(t)({\mathbf h}(t)-\mathbf e_*) =
\frac {t}{t-1}(\mathbf u - \mathbf 1)h(t)$$
and so
$${\mathbf h}(t)-\mathbf e_* =
\frac {t}{t-1}(\mathbf u - \zeta (t)^{-1}\mathbf 1)h(t).$$

Then
$$\mathbf 1^T({\mathbf h}(t)-\mathbf e_*) =
\frac {t}{t-1}(\mathbf 1^T\mathbf u - \mathbf 1^T\zeta (t)^{-1}\mathbf 1)h(t)$$
or
$$h(t)-1=
\frac {t}{t-1}(1 - \mathbf 1^T\zeta (t)^{-1}\mathbf 1)h(t).$$

Consequently, we have

\begin{lem}
$$ \frac {1-t}{h(t)}=1 - t\mathbf 1^T\zeta (t)^{-1}\mathbf 1.$$
\end{lem}

Now $N(t) = \zeta(t) - I$ is a strictly upper triangular matrix and so
$\zeta(t)$ is invertible.  In fact, $\zeta(t)^{-1} = I - N(t) + N(t)^2 -
\dots  $ and so the $(v,w)$-entry of $\zeta(t)^{-1}$ is
$$\sum_{v = v_1 > \dots > v_l = w\geq *}(-1)^{l+1}t^{|v|-|w|}.$$

Combining this remark with Lemma 1 we obtain the following result.

\begin{thm}
Let $\Gamma$ be a layered graph with unique minimal element $*$ of level $0$ and $h(t)$ 
denote the Hilbert series of $A(\Gamma)$.  Then
$$\frac{1-t}{h(t)} 
= 1  + \sum_{v_1 > v_2 \dots > v_{\ell}\ge *} (-1)^{\ell}  
t^{|v_1|-|v_{\ell}|+1}.$$
\end{thm}

We remark that the matrices $\zeta(1)$ and $\zeta(1)^{-1}$ are well-known as the zeta-matrix and 
the M\"obius-matrix of $V$ (cf. \cite R).

In the remaining sections of this paper we will use Theorem 2 to compute the Hilbert series of 
the algebras $A(\Gamma)$ associated with certain layered graphs.

\section{The Hilbert series of the algebra associated with the Hasse graph of the 
lattice of subsets of $\{1,\dots,n\}$}

Let $\Gamma_n$ denote the Hasse graph of the lattice of all subsets of $\{1,\dots,n\}$.  
Thus the vertices of $\Gamma_n$ are subsets of $\{1,\dots,n\}$, the order relation $>$ is set inclusion $\supset$, the level $|v|$ of a set $v$ is its cardinality, and the unique minimal vertex $*$ is the empty set $\emptyset$.  
Then the algebra $A(\Gamma_n)$ is the agebra $Q_n$ defined in \cite{GRW}.  
In this section we will prove the following theorem (from \cite{GGRSW}). 
The present proof is much shorter and more conceptual than that in \cite{GGRSW}.

\begin{thm}  $$H(Q_n,t) =  \frac{1-t}{1 - t(2-t)^n}.$$

\end{thm}

Our computations depend on the following lemma and corollary.
\\
\begin{lem} Let $w$ be a finite set. Then 
$$\sum_{\stackrel{w \supset w_2 \supset \dots \supset w_{\ell} = 
\emptyset}{}} (-1)^{\ell}  = 
(-1)^{|w|+1}. $$
\end{lem}

\begin{proof}  If $|w| = 1$, both sides are $+1$. Assume the result holds for all sets of cardinality $< |w|.$  Then
$$\sum_{\stackrel{w \supset w_2 \supset \dots \supset w_{\ell} = 
\emptyset}{}} (-1)^{\ell}  = 
\sum_{\stackrel{w \supset w_2 \supseteq \emptyset}{}}  
\sum_{\stackrel{w_2 \supset \dots \supset w_{\ell} = \emptyset}{}} 
(-1)^{\ell}$$
and, by the induction assumption, this is equal to
$$ \sum_{\stackrel{w \supset w_2 \supseteq \emptyset}{}} (-1)^{|w_2|}.$$
Since $$ \sum_{\stackrel{w \supset w_2 \supseteq \emptyset}{}} (-1)^{|w_2|}=
\sum_{\stackrel{w \supseteq w_2  \supseteq \emptyset}{}} (-1)^{|w_2|} - (-1)^{|w|}
= 0 + (-1)^{|w|+1} = (-1)^{|w|+1}$$
the proof is complete.
\end{proof}

\begin{cor} Let $v \supseteq w$ be finite sets.  Then
$$\sum_{v = v_1 \supset v_2 \supset \dots \supset v_{\ell} = w} 
(-1)^{\ell} = (-1)^{|v| - |w| + 1}.$$

\end{cor}
\begin{proof} Let $w'$ denote the complement of $w$ in $v$.  Sets $u$ satisfying $v \subseteq u \subseteq w$ are in one-to-one correspondence with subsets of $w'$ via the map $u \mapsto u \cap w'$.  Thus
$$\sum_{v = v_1 \supset v_2 \supset \dots \supset v_{\ell} = w} 
(-1)^{\ell} 
= \sum_{w' = v_1' \supset \dots \supset v_{\ell}' = \emptyset} 
(-1)^{\ell}.$$
By the lemma, this is $(-1)^{|w'|+1}$, giving the result.
\end{proof}

To prove the theorem we observe that

$$\sum_{\stackrel{v_1  \supset v_2 \supset \dots \supset v_{\ell} 
\supseteq v_{\ell} \supseteq  \emptyset}{\ell \geq 1}} 
(-1)^{\ell}   t^{|v_1|-|v_{\ell}|+1}= \sum_{\{1,\dots , n\} 
\supseteq v_1 \supseteq \emptyset} t^{|v_1| - |v_{\ell}| + 1} 
\sum_{v_1 \supset \dots \supset v_{\ell} \supseteq \emptyset} 
(-1)^{\ell}.$$
By  Corollary 1, this is 
$$\sum_{\{1,\dots , n\} \supseteq v_1 \supseteq v_{\ell} \supseteq 
\emptyset} t^{|v_1| - |v_{\ell}| + 1}(-1)^{|v_1| - |v_{\ell}| + 1}.$$
Let $u$ denote the compliment of $v_{\ell}$ in $v_1$ and $u'$ 
denote the complement of $u$ in $\{1,\dots,n\}.$  
Then the coefficient of $t^{k+1}$ in the above expression 
is the number of ways of choosing a $k$-element 
subset $u \subseteq \{1,\dots,n\}$ times the number of 
ways of choosing a subset $v \subseteq u'$.  
This is $\binom nk 2^{n-k}.$ Thus
$$\sum_{\stackrel{v_1  \supset v_2 \supset \dots \supset v_{\ell} 
\supseteq  \emptyset}{\ell \geq 1}} 
(-1)^l   t^{|v_1|-|v_{\ell}|+1} = \sum_{k=0}^n \binom nk 2^{n-k}(-t)^{k+1} 
= -t(2-t)^k.$$  In view of Theorem 2, this completes the proof of the 
theorem.

\section{The Hilbert series of  algebras associated with the Hasse 
graph of the lattice of subspaces of a finite-dimensional
vector space over a finite field}

We will denote by $\mathbf{L(n,q)}$ the Hasse graph of the lattice of 
subspaces of an $n$-dimensional space over the field $\mathbf{F_q}$ 
of $q$ elements.  Thus the vertices of $\mathbf{L(n,q)}$ are subspaces 
of ${\mathbf{F_q}}^n$, the order relation  $>$ is inclusion of 
subspaces $\supset$, the level $|U|$ of a subspace $U$ is its dimension,  and the unique minimal vertex $*$ is the zero subspace $(0)$.

\begin{thm}  
$$
\frac{1 - t}{H(A(\mathbf{L(n,q)}),t)} = 1 - t \sum_{\stackrel{m = 0}{}}^n  {\binom nm}_q  
(1 - t)(1 - tq) \dots (1 - tq^{n - m - 1}).$$
\end{thm}

Our proof depends on the following lemma and corollary.

\begin{lem} Let $U$ be a finite-dimensional vector space over $\mathbf{F_q}$. 
Then 
$$\sum_{\stackrel{U = U_1 \supset  U_2 \supset  \dots \supset U_{\ell} = 
(0)}{l \geq 1}} (-1)^{\ell}  = (-1)^{|U|+1} q^{\binom {|U|}2 }.$$ 
\end{lem}

\begin{proof} If  $|U| = 0$, the sum occuring in the lemma has a single term corresponding 
to $l=1,U = U_1 = (0).$  Then both sides of the expression in the lemma 
are equal to $-1$.  Now let $U$ be a finite-dimensional vector space and 
assume the result holds for all spaces of dimension less than $|U|.$  Then

$$\sum_{\stackrel{U = U_1 \supset  U_2 \supset \dots \supset U_{\ell} = 
(0)}{l \geq 1}} (-1)^{\ell} 
= \sum_{\stackrel{U = U_1 \supset  U_2}{}} \;\;  
\sum_{\stackrel{  U_2 \supset  \dots \supset U_{\ell} = (0)}
{l \geq 1}} (-1)^{\ell}.$$ 
By the induction assumption, this is equal to 
$$= \sum_{\stackrel{U \supset  U_2 }{}}  (-1)^{|U_2|} q^{\binom 
{|U_2|}2}.$$
It is well-known that the number of $m$-dimensional subspaces of the space 
$U$ is given by the $q$-binomial coefficient ${\binom {|U|}m}_q$.

Hence $$\sum_{\stackrel{U = U_1 \supset  U_2 \supset  \dots \supset 
U_{\ell} = (0)}{l \geq 1}} (-1)^{\ell}  
= \sum_{\stackrel{|U_2| = 0 }{}}^{|U|-1}  {\binom {|U|}{|U_2|}}_q  
(-1)^{|U_2|}q^{\binom {|U_2|}2}.$$

Recall the $q$-binomial theorem 

$$ \prod_{i = 0}^{m - 1}  (1 + xq^i) = \sum_{\stackrel{j = 0 }{}}^{m}  
{\binom mj}_q   q^{\binom j2} x^j.$$

Set $x = -1$. Then the $i = 0$ factor in the product is $0$ and so we have

$$\sum_{\stackrel{j = 0 }{}}^{m - 1}  {\binom mj}_q  (-1)^j  
q^{\binom j2} = (-1)^{m+1} q^{\binom m2}.$$

Thus

$$\sum_{\stackrel{U = U_1 \supset  U_2 \supset \dots \supset U_{\ell} = 
(0)}{l \geq 1}} (-1)^{\ell}  = (-1)^{|U| + 1} q^{\binom {|U|}2}$$
as required.
\end{proof}

\begin{cor} Let $V \supseteq W$ be subspaces of ${\mathbf F_q}$. Then
$$\sum_{V = V_1 \supset V_2 \supset \dots \supseteq V_{\ell} =W} 
(-1)^{\ell} = 
(-1)^{|V/W|+1}q^{{\binom {|V/W|}2}}.$$

\end{cor}
\begin{proof}
Since subspaces $Y, V \supseteq Y \supseteq W$, are in one-to-one correspondence with subspaces of $V/W$ 
via the map $Y \mapsto Y/W$, this is immediate from the lemma.
\end{proof}

To prove the theorem, we observe that 
$$\sum_{\stackrel{V_1 \supset V_2 \supset \dots \supset V_{\ell} \supseteq 
(0)}{l \geq 1}} (-1)^{\ell}t^{|V_1/V_{\ell}|+1} =
\sum_{{\mathbf{F_q}}^n \supseteq V_1 \supseteq V_{\ell} \supseteq (0)} 
t^{|V_1/V_{\ell}|+1} \sum_{\stackrel{V_1 \supset V_2 \supset \dots \supset 
V_{\ell} \supseteq (0)}{\ell \geq 1}} (-1)^{\ell}.$$
By Corollary 2, this is equal to
$$\sum_{{\mathbf{F_q}}^n \supseteq V_1 \supseteq V_{\ell} \supseteq (0)} 
t^{|V_1/V_{\ell}|+1}(-1)^{|V_1/V_{\ell}|+1}q^{{\binom{|V_1/V_{\ell}|}2}}.$$

Set $|v_{\ell}| = m$ and $|V_1/V_{\ell}| = k$.  Then the number of 
possible $V_{\ell}$ is $\binom nm_q$ and, for fixed $V_{\ell}l$, the 
number 
of 
possible $V_1$ is the number of $k$-dimensional subspaces of 
${\mathbf{F_q}}^n/V_{\ell}$ which is $\binom {n-m}k_q.$  Thus

$$\sum_{\stackrel{V_1 \supset V_2 \supset \dots \supset V_{\ell} \supseteq 
(0)}{\ell \geq 1}} (-1)^{\ell}t^{|V_1/V_{\ell}|+1} =
 \sum_{\stackrel{0 < k,m}{k+ m 
\leq 
n}}  {\binom nm}_q  {\binom {n-m}k}_q  (-t)^{k+1} q^{\binom k2}$$

$$ = (-t)\sum_{m=0}^n  
{\binom nm}_q \sum_{k=0}^{n-m} {\binom {n - m}k}_q  (-t)^k q^{\binom k2}.$$  

Setting $x = -t$ in the $q$-binomial theorem shows that 
$$\sum_{k=0}^{n-m} {\binom {n - m}k}_q  (-t)^k q^{\binom k2} = \prod_{i=0}^{n-m-1} (1 - tq^i).$$
Therefore 
$$\sum_{\stackrel{V_1 \supset V_2 \supset \dots \supset V_{\ell} \supseteq 
(0)}{\ell \geq 1}} (-1)^{\ell}t^{|V_1/V_{\ell}|+1} =
(-t)\sum_{m=0}^n \binom nm_q \prod_{i=0}^{n-m-1} (1 - tq^i).$$
In view of Theorem 2, the theorem is proved.

Note that setting $q = 1$ in the expression in Theorem 4 gives $ 1 - 
t(2-t)^n.$  By Theorem 3, this is $\frac{1-t}{H(Q_n,t)}.$

Recall (cf. \cite{U}) that if $A$ is a quadratic algebra it has a dual quadratic algebra, denoted $A^!$ 
and that if $A$ is a Koszul algebra 
the Hilbert series of $A$ and $A^!$ are related by 
$$H(A,t)H(A^!,-t) = 1$$ 

Since by \cite{RSW}  $A(\mathbf{L(n,q)})$ is a Koszul algebra, we have the following

\begin{cor} $$H(A(\mathbf{L(n,q)})^!,t) = 1 + \sum_{m=0}^{n-1} {\binom nm_q }(1 + tq)\dots(1+tq^{n-m-1}).$$
\end{cor}

\section{The Hilbert series of algebras associated with complete layered graphs}

We say that a layered graph $\Gamma= (V,E)$ with $V = \cup _{i=0}^nV_i$ 
is {\bf complete} if for every $i,
1 \le i \le n,$ and every $v \in V_1, w \in V_{i-1}$ , there is a unique edge $e$ 
with ${\mathbf t}(e) = v, {\mathbf h}(e) = w.$  
A complete layered graph is determined
(up to isomorphism) by the cardinalities of the $V_i$.  
We denote the complete layered graph 
with $V = \cup_{i=0}^n V_i, |V_i| = m_i$ for $0 \le i \le n$, by
${\mathbf{C}}[m_n,m_{n-1},\dots ,m_1,m_0]$.  Note that the graph 
${\mathbf{C}}[m_n,m_{n-1},\dots ,m_1,1]$
has a unique minimal vertex of level $0$ and so Theorem 2
applies to 
$A({\mathbf{C}}[m_n,m_{n-1},\dots ,m_1,1])$.  We will show:

\begin{thm}  $$\frac{1-t}{H(A({\mathbf{C}}[m_n,m_{n-1},\dots ,m_1,1],t)} = $$
$$ 1 - \sum_{k = 0}^n  \sum_{a = k}^n  \; (-1)^k 
m_a (m_{a-1} - 1)(m_{a-2} - 1) \dots (m_{a-k+1} - 1)m_{a-k} t^{k+1}.$$
\end{thm}

\begin{proof}

We first compute   

$$\sum_{\stackrel{v_1  > v_2 > \dots > 
v_{\ell} \geq  \emptyset}{\ell \geq 1}} (-1)^{\ell}   
t^{|v_1|-|v_{\ell}|+1}.$$  

The coefficient of $t^{k+1}$ in the sum is 
  
$$\sum_{\stackrel{v_1 > v_2 > \dots > v_{\ell} \geq *}
{\ell \geq 1, \; |v_1| - |v_{\ell}| = k}} (-1)^{\ell}  
= \sum_{|v_1|=k}^n \sum_{\stackrel{v_1 > \dots > 
v_{\ell}}{|v_1|-|v_{\ell}| = k}
} (-1)^{\ell}.$$

Note that the number of chains $v_1 > \dots > v_{\ell}$ with $|v_i| = a_i$ 
for $1 \le i \le \ell$ is $m_{a_1}m_{a_2}\dots m_{a_{\ell}}.$  Then, 
writing 
$k = |v_1| - |v_{\ell}| $ and $a_1 = a$ we have 
$$\sum_{\stackrel{v_1  > v_2 > \dots > 
v_{\ell} \geq  *}{\ell \geq 1}} (-1)^{\ell}   t^{|v_1|-|v_{\ell}|+1} =
\sum_{k=0}^n \left ( \sum_{\stackrel{v_1 > \dots > v_{\ell} \ge *}{\ell 
\ge 1, \; |v_1| - |v_{\ell}|=k}} (-1)^{\ell} \right ) t^{k+1}$$
$$ = \sum_{k=0}^n\left (\sum_{\stackrel{a_1> \dots a_{\ell-1} > a_1 - k 
\ge 0}{\ell \ge 1} }
(-1)^l m_{a_1} m_{a_2} \dots m_{a_{\ell-1}} m_{a_{1 - k}}\right )t^{k+1}$$
$$ = \sum_{k=0}^n\left (\sum_{a = k}^n m_a(1 - m_{a-1}) \dots (1 - m_{a-k+1})m_{a-k}\right )t^{k+1}. 
$$

The theorem now follows from Theorem 2.
\end{proof}
This result applies, in particular, to the case $m_0 = m_1 = ... = m_n = 1.$  
The resulting algebra $A({\mathbf C}[1,\dots,1])$ has $n$ generators and no relations.  
Theorem 5 shows that
$$\frac{1-t}{H(A({\mathbf C}[1,\dots,1]),t)} 
= 1 - \sum_{a=0}^n t + \sum_{a=1}^n t^2 = (1-t)(1-nt).$$
Thus $H(A({\mathbf C}[1,\dots,1]),t) = \frac1{1-nt}$ and 
we have recovered the well-known expression for the Hilbert 
series of the free associative algebra on $n$ generators.

Since by \cite{RSW}  the algebras associated to complete directed graphs 
are Koszul algebras,  we have the following 
\\

\begin{cor}
$$H(A({\mathbf{C}}[m_n,m_{n-1},\dots ,m_1,1])^!,t)= 
$$ $$1 +     \sum_{k=1}^n\sum_{a = 
k}^n  \;  
m_a (m_{a-1} - 1)(m_{a-2} - 1) \dots (m_{a-k+1} - 1)\; t^k.$$
\end{cor}

\end{document}